\catcode`\@=11

\font\tenmsa=msam10 \font\sevenmsa=msam7 \font\fivemsa=msam5
\font\tenmsb=msbm10 \font\sevenmsb=msbm7 \font\fivemsb=msbm5
\newfam\msafam \newfam\msbfam
\textfont\msafam=\tenmsa  \scriptfont\msafam=\sevenmsa
\scriptscriptfont\msafam=\fivemsa
\textfont\msbfam=\tenmsb  \scriptfont\msbfam=\sevenmsb
\scriptscriptfont\msbfam=\fivemsb

\font\teneufm=eufm10 \font\seveneufm=eufm7 \font\fiveeufm=eufm5
\newfam\eufmfam
\textfont\eufmfam=\teneufm
\scriptfont\eufmfam=\seveneufm\scriptscriptfont\eufmfam=\fiveeufm

\def\hexnumber@#1{\ifcase#1 0\or1\or2\or3\or4\or5\or6\or7\or8\or9\or
 A\or B\or C\or D\or E\or F\fi }

\edef\msa@{\hexnumber@\msafam} \edef\msb@{\hexnumber@\msbfam}

\mathchardef\square="0\msa@03

\def\Bbb#1{{\fam\msbfam\relax#1}}

\catcode`\@=12



\magnification\magstephalf\nopagenumbers
\vsize=19truecm \hsize=14truecm \voffset=1.5cm \hoffset.8cm
\looseness=1 \tolerance=500

\font\caps=cmcsc10


\font\title=cmbx10 scaled \magstep 1 
\font\section=cmr10 scaled \magstep 2
\font\eighttimesit=cmr8

\def\makeheadline{\vbox to 0pt{\vskip-25pt \line{\vbox
to8.5pt{}\the\headline}\vss}\nointerlineskip}
\def\makefootline{\vbox to 0pt{\vskip1pt \line{\vbox
to8.5pt{}\the\footline}\vss}\nointerlineskip}
\def\endsection{\vskip1cm\goodbreak}

\parskip=5pt


\def\P{{\Bbb P}}  
\def\qed{\hfill\vbox{\hrule\hbox{\vrule\kern3pt
    \vbox{\kern6pt}\kern3pt\vrule}\hrule}} 

\def\extendibilidad{1.1}
\def\producto{1.2}
\def\degenerada{1.3}
\def\Veronese{1.4}
\def\unodos{1.5}
\def\dosdos{1.6}


\def\proydual{2.1}
\def\muchasrectas{2.2}
\def\clas-proy{2.3}
\def\fibrado-desc{2.4}

\def\Yo{1}
\def\ABTcurvafundamental{2}
\def\AS{3}
\def\Rogora{4}
\def\Segre{5}
\def\Severi{6}
\def\Zak{7}

\centerline{\title The universal rank-$(n-1)$ bundle
on $G(1,n)$ restricted to subvarieties}
\bigskip
\centerline{{\caps Enrique Arrondo}\footnote{(*)}{This work was supported in
part by  DGICYT grant No. PB93-0440-C03-01}}
\bigskip

\centerline{\eighttimesit Departamento de Algebra}
\centerline{\eighttimesit Facultad de Matem\'aticas, Universidad Complutense,
28040 Madrid, Spain}
\centerline{\eighttimesit E-mail: enrique@sunal1.mat.ucm.es}

\bigskip
\bigskip

\centerline{\sl En record de Ferran Serrano, 
extraordinari com a persona i com a matem\`atic}

\bigskip
\bigskip

\centerline{\caps Abstract} 

\hsize 12truecm
\itemitem{}We classify those smooth $(n-1)$-folds in $G(1,\P^n)$ for which the
restriction of the rank-$(n-1)$ universal bundle has more than $n+1$
independent sections. As an aplication, we classify also those $(n-1)$-folds
for which that bundle splits.

\hsize 14truecm

\bigskip
\bigskip

\noindent It is a classical problem to study which projective varieties of
small codimension are not linearly normal (i.e. isomorphically projected from
higher projective spaces). The first observation is that any $n$-dimensional
variety can be projected to $\P^{2n+1}$, but is expected to produce singular
points when projected to $\P^{2n}$. Hence, $n$-dimensional varieties of
codimension at most $n$ are expected to be linearly normal. For $n=2$, Severi
proved (see his classical paper [\Severi]) that the Veronese surface is the
only smooth surface in $\P^5$ that can be isomorphically projected to
$\P^4$. In general, projectability (or linear normality) is characterized by
the dimension of the secant varieties. A recent thorough study of secant
varieties can be found in [\Zak], where a lot of projectability results are
given. 

More generally, one can study which varieties are isomorphically mapped under
the projection from a Grassmannian $G(k,N)$ of $k$-planes in $\P^n$ to another
$G(k,n)$ induced by linear projections from $\P^N$ to $\P^n$. In the
particular case $k=1$, we have that $G(1,n)$ has even dimension $2(n-1)$. Thus,
it is natural to ask which $(n-1)$-dimensional varieties in $G(1,n)$ are
projected from a bigger $G(1,N)$. This problem was solved in [\AS]\ for $n=3$,
and then studied in [\Yo]\ for any dimension. In the latter, we introduced an
appropriate notion of secant varieties for Grassmannians of lines. However,
the theory only worked if we added the hypothesis of ``uncompressedness'',
i.e. that the union of the lines parametrized by our varieties has the
expected dimension. 

The situation in the other extreme case, when $k$ is very big, is just the
opposite, since the union of the corresponding $k$-planes can never have the
expected dimension, except for varieties of big codimension). The
motivation of this paper was to understand this new phenomenon of
compressedness that makes the theory so different for a general
$k$. To this purpose, we will study the maximal values of $k$. In
particular, we will study which varieties of $G(k,n)$ are projected from
higher Grassmannians for $k=n-1,n-2$. The case $k=n-1$ is easy (see
Proposition \proydual), so that we will study more in detail the case
$k=n-2$. Since $G(n-2,n)$ has again even dimension $2(n-2)$, we will also
specially study its $(n-1)$-dimensional subvarieties that come from a bigger
$G(n-2,N)$. Another reason to study this case is that to a variety in
$G(n-2,n)$ corresponds by duality another variety in $G(1,n)$. In fact we will
write our results in terms of varieties in $G(1,n)$. Varieties of dimension
$n-1$ in $G(1,n)$ (particularly the cases $n=3,4$) have been
thoroughly studied by the classical geometers one century ago, and this
research has retaken interest nowadays. It should also be mentioned that 
projectability of the dual varieties can also be interpreted as some kind
extendability of the original variety (see Remark \extendibilidad).

To my knowledge, for $n\ge 4$, few things are known about the dual varieties
of varieties in $G(1,n)$ (or more likely, the classical geometers knew but
never wrote up). This is why we needed to devote most of the first section to
not only recall classical examples of those $(n-1)$-folds, but also to
carefully describe their dual varieties in some case. 

In section 2, we give a characterization (see Theorem \clas-proy) of which
$(n-1)$-folds of $G(1,n)$ verify that its dual variety in $G(n-2,n)$ is
projected from a higher Grassmannian of $(n-2)$-planes. On the other hand, since projectability is related to the restriction
of the rank-$(n-1)$ universal bundle, we are able to characterize when such a
restriction to an $(n-1)$-fold splits (see Proposition \fibrado-desc). The same
problem was solved in [\AS] for the restriction of the rank-two universal
bundles of $G(1,3)$ to surfaces (in fact we took our method from there). Such
a classification can be considered as the very first step towards the
difficult problem of studying the stability of the restriction of the
universal bundles of Grassmannians to subvarieties. 

In order to deal with the intermediate values of $k$, it would be needed to
have a full knowledge of varieties containing more linear spaces than
expected. As far as I know, these kind of results are known only when the
dimension of the linear spaces is close to the maximum allowed. There are
results by Rogora (see [\Rogora]) for varieties with a lot of lines, and
for general linear spaces, there are classical results by B. Segre (see
[\Segre], or Lemma\muchasrectas for a particular case), which we used for the
case $k=n-1$. I hope that a development of such a theory will help in the
future to deal with the projectability theory for any value of $k$.

\endsection

\centerline{\section 1. Preliminaries and examples}
\bigskip

\noindent{\sl Notation.} All the varieties we consider will be complex and
integral and, unless otherwise specified, they are also assumed to be smooth.
For a vector bundle $E$ over a variety $X$, $H^i(X,E)$ means the space of global sections of a vector bundle
over a variety $X$, while $h^i(X,E)$ will denote the dimension of that space.
We denote by $G(k,n)$ the Grassmannian of $k$-linear spaces in the complex
projective space $\P^n$. More generally, if
$\Lambda$ is a linear subspace of $\P^n$, $G(k,\Lambda)$ will denote the set
of $k$-linear subspaces of $\Lambda$. We denote by ${\cal Q}$ and ${\cal S}$
the universal vector bundles of respective ranks $k+1$ and $n-k$. More
precisely, ${\cal Q}$ and ${\cal S}$ will be the vector bundles appearing in
the universal exact sequence
$$
0\to\check{\cal S}\to H^0(\P^n,{\cal O}_{\P^n}(1))\to {\cal Q}\to 0
$$
This means that, if $Y$ is a subvariety of $G(k,n)$, then the restriction
${\cal Q}_{|Y}$ is the vector bundle that embeds $Y$ in $G(k,n)$, while ${\cal
S}_{|Y}$ embeds $Y$ in the dual Grassmannian $G(n-k-1,n)=G(n-k-1,{\P^n}^*)$.
The image of $Y$ in $G(n-k-1,{\P^n}^*)$ will be denoted by $Y^*$, and will be
called the {\it dual subvariety} of $Y$. For us, a {\it congruence} in
$G(1,n)$ will mean a smooth irreducible subvariety of $G(1,n)$ of dimension
$n-1$. By abuse of notation, a $k$-plane of a variety $Y\subset G(k,n)$ will
just mean a $k$-plane in $\P^n$ corresponding to a point of $Y$. 
\bigskip

\noindent{\sl Remark \extendibilidad.} Assume that a subvariety $Y\subset
G(k,n)$ verifies that $h^0(Y,{\cal S}_{|Y})\ge n+2$. This means that the dual
$Y^*\subset G(n-k-1,n)$ is an isomorphic projection of a subvariety
$Y'\subset G(n-k-1,n+1)$ (this projection being induced by a linear projection
from ${\P^{n+1}}^*$ to ${\P^n}^*$) which is not contained in $G(n-k-1,H)$ for
any hyperplane $H\subset{\P^{n+1}}^*$. In this situation, we will say that
$Y^*$ is a {\it nontrivial projection} of $Y'$, or that $Y^*$ {\it is
projected} from $G(n-k-1,n+1)$. Dually, this means that there exists a
subvariety $\hat Y\subset G(k+1,n+1)$, parametrizing $(k+1)$-planes in
$\P^{n+1}$ not all of them passing through a point, such that the intersection
of those $(k+1)$-planes with a hyperplane $\P^n$ produces exactly our family
$Y$ of $k$-planes. Therefore, in this sense, $H^0(Y,{\cal S}_{|Y})$ measures
how far $Y$ can be extended.

\bigskip

We look now at different examples of congruences in $G(1,n)$. These will be
the examples appearing in the statements of our main results in the next
section.

\bigskip

\noindent{\bf Example \producto.} Let us consider on $Y=\P^r\times\P^{n-1-r}$
the vector bundle ${\cal O}_Y(1,0)\oplus{\cal O}_Y(0,1)$. This defines an
embedding of $Y$ as a congruence in $G(1,n)$ consisting of the set of lines
joining two disjoint linear spaces of dimensions $r$ and $n-1-r$. The dual of
${\cal S}_Y$ is the kernel of the epimorphism
$$
H^0(Y,{\cal O}_Y(1))\otimes{\cal O}_Y\to {\cal O}_Y(1,0)\otimes{\cal O}_Y(0,1)
$$
Since $H^0(Y,{\cal O}_Y(1))$ splits as $H^0(\P^r,{\cal O}_{\P^r}(1))\oplus
H^0(\P^{n-r-1},{\cal O}_{\P^{n-r-1}}(1))$, it immediately follows that
${\cal S}_Y\cong pr_1^*(T_{\P^r}(-1))\oplus pr_2^*(T_{\P^{n-r-1}}(-1))$.
\bigskip

\noindent{\bf Example \degenerada.} Assume that $Y\subset G(1,n)$  is
contained in a $G'=G(1,H)$ for some hyperplane $H\subset\P^n$. If ${\cal Q}'$
and ${\cal S}'$ are the universal bundles on $G'$, it holds that ${\cal
Q}_{|G'}={\cal Q}'$ and ${\cal S}_{|G'}\cong {\cal S}'\oplus{\cal O}_{G'}$. In
particular, ${\cal S}_{|Y}\cong {\cal S}'_{|Y}\oplus{\cal O}_Y$. Reciprocally,
if ${\cal S}_{|Y}\cong E\oplus{\cal O}_Y$ for some rank-$(n-2)$ vector bundle
$E$, then $Y$ is contained in some $G(1,H)$. Indeed, such a decomposition
implies that all the $(n-2)$-planes of the dual congruence $Y^*$ pass through
a fix point. Dualizing, all the lines of $Y$ are contained in a hyperplane $H$.

\bigskip

\noindent{\bf Example \Veronese.} Let $Q$ be a smooth quadric in $\P^5$. It
is well-known that $Q$ has two families of planes, each of them parametrized
by $\P^3$. The map from $\P^3$ to $G(2,\P^5)$ defining any of the two
families is given by the vector bundle $\Omega_{\P^3}(2)$. The dual
Grassmannian is another $G(2,5)$ and it also holds that the dual of the set of
planes in $Q$ is again the set of planes in the dual quadric
$Q^*\subset{\P^5}^*$. If we take a general hyperplane $\P^4$ in $\P^5$, then we
obtain a smooth three-dimensional quadric $Q'$, and any line in $Q'$ is the
intersection of $\P^4$ with a unique plane in $Q$ in each family. By Remark
\extendibilidad, the dual congruence $Y^*\subset G(2,4)$ of the congruence
$Y\subset G(1,4)$ of lines in $Q'$ is projected from a threefold
$Y'\subset G(2,5)$, and this $Y'$ is dual to one of the two family of planes in
$Q$. Hence, $Y'$ itself consists of the set of planes contained in a smooth
quadric. 

\bigskip

Finally we recall now two examples from [\ABTcurvafundamental].
\bigskip

\noindent{\bf Example \unodos.} Let us consider the congruence $Y$ of
$G(1,n)$ studied in [\ABTcurvafundamental] 3.3. Geometrically, it is described
in the following way. In $\P^n$ there are a plane $\Pi'$ and a linear space
$\Lambda$ of dimension $n-2$ meeting in one point $P$ and there is a smooth
conic $C\subset\Pi'$ passing through $P$. Then the congruence consists of the
closure of the set of lines joining a point of $\Lambda\setminus\{P\}$ and a
point of $C\setminus\{P\}$. In order to describe its dual congruence, we look
at the following alternative description. We choose in $C$ a point $P_0$
different from $P$. The stereographic projections of $C$ from $P$ and $P_0$
define an isomorphism $\varphi$ between the pencil $A$ of hyperplanes
containing $\Lambda$ and the pencil of lines in $\Pi'$ passing through $P_0$.
Specifically, a hyperplane in the pencil $A$ meets $C$ in $P$ plus another
point $Q$, and we define the image of that hyperplane to be the line $P_0Q$ in
the second pencil. The congruence is then a scroll parametrized by $A$ in the
following way. For any hyperplane $H\in A$, the lines of the congruence
associated with it are those obtained intersecting $H$ with the set of planes
containing the line $\varphi(H)$.

This immediately leads to the following description of the dual congruence.
There is a line $L\subset{\P^n}^*$ and an isomorhism $\varphi$ from $L$ to
the pencil of $(n-2)$-planes contained in a hyperplane $H'\subset{\P^n}^*$ and
containing an $(n-3)$-plane $\Lambda'\subset H'$. Then, for any point $P\in L$,
the $(n-2)$-planes of $Y^*$ associated with it are those spanned by $P$ and the
set of $(n-3)$-planes contained in $\varphi(P)$. Our next task is to describe
this in terms of an embedding induced by a rank-$(n-1)$ vector bundle on an
$(n-1)$-fold.

A pencil of $(n-2)$-planes is a map from $\P^1$ to $G(n-2,n)$ induced by
the vector bundle ${\cal O}_{\P^1}^{\oplus n-2}\oplus{\cal O}_{\P^1}(1)$.
Therefore, the set of $(n-3)$-planes in some $(n-2)$-plane of the above pencil
is parametrized by the projective bundle $X=\P({\cal O}_{\P^1}^{\oplus
n-2}\oplus{\cal O}_{\P^1}(-1))$ (i.e. the relative dual projective bundle to
the previous one). Let $p:X\to\P^1$ be the structure morphism and denote by
${\cal O}_X(1)$ its tautological line bundle. Then the embedding of $X$ in
$G(n-3,n)$ that associates to each point of $X$ the $(n-3)$-plane in $\P^n$ it
represents is given by the vector bundle $T_{X/\P^1}(-1)$. Hence, $Y^*$ is the
image in $G(n-2,n)$ of a morphism induced by the rank-$(n-1)$ vector bundle
$T_{X/\P^1}(-1)\oplus p^*{\cal O}_{\P^1}(1)$. Observe that this vector bundle
has $n+2$ independent sections, so that it embeds $X$ in $G(n-2,n+1)$, and
$Y^*$ is just a linear projection of it. In fact, the description in
$G(n-2,n+1)$ coincides with the one given for $Y^*$, but now the hyperplane
$H'$ and the line $L$ are disjoint.

\bigskip

\noindent{\bf Example \dosdos.} Finally, we recall know from 
[\ABTcurvafundamental] 3.2 an example of another congruence $Y\subset G(1,n)$
that will be of interest for us. We will not describe it here completely as
in the previous example, but we will just pay attention to the few facts that
will be needed. For this congruence, there is a conic $C\subset\P^n$ such
that, for each point $P\in C$, the set of lines of $Y$ passing through it is
the set of lines contained in a hyperplane $H_P\ni P$ and passing through it.
Then, $Y$ has a structure of a scroll over $C\cong\P^1$, and, under the
Pl\"ucker embedding, it is in fact a rational normal scroll of degree $2n-2$.

One way of constructing a congruence of this type is by taking in $G(1,n)$ the
dependency locus of three sections of the vector bundle ${\cal
Q}\oplus{\cal S}$. It is not difficult to check --by looking at its
invariants-- that this dependency locus is indeed a congruence of the desired
type. However, a more intuitive way of seeing the scroll structure in this
construction is the following. A section of ${\cal Q}\oplus{\cal S}$ vanishes
on the set of lines contained in a hyperplane $H\subset\P^n$ and passing
through a point $P\in\P^n$. If the section is general, $P\notin H$ and hence
this locus is empty; but for a special section, $P$ belongs to $H$ and the
zero locus of the section is a Schubert variety isomorphic to
$\P^{n-2}$. If we take a general net of sections of ${\cal Q}\oplus{\cal S}$,
we have a conic inside that net for which the corresponding section has a
nonempty zero locus. This gives the scroll structure over a conic. From that
construction, it is easy to see that $h^0(Y,{\cal S}_{|Y})=n+2$, which implies
that $Y^*$ is projected from $G(n-2,n+1)$.

\endsection

\centerline{\section 2. Projection of subvarieties of Grassmannians}

\bigskip

\noindent We first study the easiest Grassmannians, i.e. those which are in
fact a dual projective space. We will see that the fact that a subvariety of
any dimension in $G(n-1,n)$ there is projected from a higher Grassmannian is
impossible except for curves, and will study in that case what is the
situation for congruences (i.e. $n=2$). More precisely, we have the following
\bigskip

\noindent {\bf Proposition \proydual} 

{\sl If a subvariety $Y\subset
G(n-1,n)$ is projected from $G(n-1,n+1)$, then $Y$ is a curve. Morevoer, if
$n=2$, then $Y$ is a conic. Precisely, $Y$ consists of the set of tangent lines
to a conic, and it is projected from the subvariety of $G(1,3)$ parametrizing
the set of lines of one of the rulings of a smooth quadric in $\P^3$; this 
is not projected from $G(1,4)$.}
\bigskip

\noindent{\sl Proof.}
Let $Y$ be a subvariety of $G(n-1,n)$ that is projected from a subvariety 
of $G(n-1,n+1)$. Writing $\P^n= G(n-1,n)$, this means that
$h^0(Y,T_{\P^n}(-1)_{|Y})>n+1$. But from the dual of the Euler exact sequence
restricted to $Y$,
$$
0\to{\cal O}_Y(-1)\to H^0(\P^n,{\cal O}_{\P^n}(1))^*\otimes{\cal O}_Ys \to
T_{\P^n}(-1)_{|Y}\to 0,
$$
it follows that $h^1(Y,{\cal O}_Y(-1))\neq 0$. Then, by the Kodaira vanishing
theorem, $Y$ must be a curve. If, moreover, $n=2$, $Y$ is a
plane curve of some degree $d$, and from the exact sequence 
$$
0\to T_{\P^2}(-1-d)\to T_{\P^2}(-1)
\to T_{\P^2}(-1)_{|Y}\to 0$$
it follows that $h^0(Y,T_{\P^2}(-1)_{|Y})>3$ if and only if
$h^1(Y,T_{\P^2}(-1-d))\neq 0$. But it is clear that the latter holds if
and only if $d=2$. Also, in this case 
$h^0(Y,T_{\P^2}(-1)_{|Y})=4$, so that $Y$ comes only from
$G(1,3)$. In fact, $T_{\P^2}(-1)_{|Y}\cong {\cal O}_{\P^1}(1)^{\oplus 2}$, so
that the corresponding curve in $G(1,3)$ parametrizes the lines in one of the
rulings of a smooth quadric. 
\qed

\bigskip

\noindent{\sl Remark.} Let us give a more geometric proof of the above
result to illustrate in part what will be our strategy for the case of
subvarieties in $G(n-2,n)$. Assume $Y\subset G(n-1,n)$ is projected from a
variety $Y'\subset G(n-1,n+1)$. In particular, the union in $\P^{n+1}$ of all
the $(n-1)$-planes parametrized by $Y'$ must be a hypersurface $X$. But, if
dim$Y\ge 2$, a variety $X$ with so many linear spaces is necessarily a
hyperplane. Indeed, through a general point of $X$ there will be at least a
one-dimensional family of $(n-1)$-planes of $Y'$; hence, through two general
points of $X$ there passes an $(n-1)$-plane contained in $X$, and this
implies that $X$ is linear. Therefore, the projection is trivial. On the other
hand, if $Y$ is a curve and $n=2$, we recall from Remark \extendibilidad\ that
this means that $Y^*$ is the plane section of a ruled surface in $\P^3$. But
it is a well-known fact that a ruled surface of degree $d$ in ${\bf P}^3$ has a
singular curve of degree $d-2$, so that its hyperplane section can be smooth
if and only if $d=2$.

\bigskip

Hence, in general we need to study projective varieties with a lot of linear
spaces. This has been done by B. Segre in [\Segre]. There is a recent paper by
Rogora (see [\Rogora]) studying projective varieties with too many lines. For
the reader's convenience, we will proof here, by using the result of Rogora,
the part we will need of the classical result of Segre. I want to thank Dario
Portelli for providing me the old paper by Segre.

\bigskip

\noindent {\bf Lemma \muchasrectas} 

{\sl Let $X$ be a (not necessarily smooth) hypersurface of $\P^{n+1}$ which
is the union of an irreducible $(n-1)$-dimensional family $Y'$ of
$(n-2)$-planes. If $n\ge4$, then either $X$ is a hyperplane, or a quadric or
the $(n-2)$-planes are distributed in a one-dimensional family of
$(n-1)$-planes.  }
\bigskip

\noindent{\sl Proof.} By Rogora's result, it suffices to show that $X$
contains at least a $(2n-3)$-dimensional family of lines. We look first at
the incidence variety
$$
I=\{(\Lambda,L)\in G(n-2,n+1)\times G(1,n+1)\ | \ \Lambda\in Y',\
L\subset\Lambda\}
$$

The projection of $I$ over $Y'$ has fibers isomorphic to $G(1,n-2)$.
Therefore, $I$ is irreducible of dimension $3n-7$. Since we want to prove that
the image of $I$ under the second projection has dimension at least
$2n-3$, it is enough to show that a general line contained in a general
$(n-2)$-plane of $Y'$ is contained in at most an $(n-4)$-dimensional family
of $(n-2)$-planes of $Y'$. 

I claim first that two general $(n-2)$-planes of $Y'$ do not meet
along an $(n-3)$-plane. Otherwise, either all the $(n-2)$-planes of $Y'$ are
in the same $(n-1)$-plane, or they all contain the same $(n-3)$-plane $A$. In
the latter case, we observe that the set of $(n-2)$-planes of $\P^{n+1}$
containing $A$ has dimension three; and since $n\ge4$, $Y'$ must consist of the
whole set of those planes. In either case, the union of the $(n-2)$-planes of
$Y'$ is not a hypersurface of ${\bf P}^{n+1}$, so that we get a contradiction.

We take now a general $(n-2)$-plane $\Lambda$ of $Y'$. Because of the
above claim, the intersection of $\Lambda$ with the rest of $(n-2)$-planes of
$Y'$ will consist only of a family $Y'_1$ of dimension $r\le n-2$ of
$(n-3)$-planes and a family $Y'_2$ of dimension $n-1$ of $(n-4)$-planes.
When we take a general line $L\subset\Lambda$, the dimension of the
$(n-3)$-planes of $Y'_1$ containing $L$ is then $r-2$, while the dimension of
the the $(n-4)$-planes of $Y'_2$ containing $L$ is $n-5$. Hence the dimension
of $(n-2)$-planes of $Y'$ containing $L$ is at most $n-4$, as wanted.
\qed
\bigskip

\noindent {\bf Theorem \clas-proy}

{\sl Assume $s\ge n-1$, $n\ge4$ and let $Y\subset G(1,n)$ be a smooth
variety of dimension $s$ such that its dual $Y^*\subset G(n-2,n)$ is a
nontrivial projection from $G(n-2,n+1)$. Then $s=n-1$ and $Y$ is one of the
following:
\item{a)} The congruence of lines in a smooth quadric of $\P^4$.
\item{b)} The congruence of Example \unodos.
\item{c)} The congruence of Example \dosdos.

\noindent Moreover, in all three cases, $Y^*$ is projected from $G(n-2,n+1)$,
but not from
$G(n-2,n+2)$}

\bigskip
\noindent{\sl Proof.} Let $Y'\subset G(n-2,n+1)$ be a smooth variety of
dimension $s$ that can be projected isomorphically to $Y^*\subset G(n-2,n)$ and
let
$X\subset{\P^{n+1}}^*$ be the union of all the $(n-2)$-planes of $Y'$. Since
$Y'$ is  projectable, this implies in particular that $X$ is not the whole
$\P^{n+1}$, so that it has dimension at most $n$ and is not a linear
space. Clearly, $X$ has dimension exactly $n$, since otherwise two general
points of it would lie in a linear space contained in $X$, and $X$ would be
linear. By Lemma \muchasrectas, either $X$ is quadric or the $(n-2)$-planes are
distributed in a one-dimensional family of $(n-1)$-planes (for the last case,
we apply the lemma for all subvarieties of $Y'$ of dimension $n-1$, and
observe that $X$, not being linear, cannot contain a two-dimensional family of
$(n-1)$-planes).

In the first case, we first observe that an $n$-dimensional quadric that is a
cone over a smooth $r$-dimensional quadric (and hence it has a vertex of
dimension $n-r-1$) contains linear spaces of dimension at most $n-{r\over
2}$. This implies that $r\le 4$. If $r=4$, then the family of $(n-2)$-planes
contained in the quadric has dimension three, so that it has to be $s=3$
and $n=4$. Hence, by dimensional reasons, $Y'$ must coincide with one of the
two families of planes in a smooth quadric of $\P^5$. This is then example
\Veronese, so that $Y$ is as stated in a). If $r=3$, the family of
$(n-2)$-planes contained in the quadric has dimension three, so that
again $s=3$ and $n=4$, and it holds now that all the planes contained in the
quadric are in a one-dimensional family of three-dimensional linear spaces.
Finally, if $r\le 2$, then $n$ can take any value, but it always happens that
the $(n-2)$-planes contained in a quadric are in a one-dimensional family of
linear spaces of dimension $n-1$. Rather that studying separately these
last two subcases, we will consider them as a particular case of the following.

So we can assume from now on that the $(n-2)$-planes of $Y'$ are distributed in
a one-dimensional family of $(n-1)$-planes contained in $X$. We first observe
that this property will be preserved under projection, i.e. the $(n-2)$-planes
of $Y^*$ are distributed in a one-dimensional family of hyperplanes of
${\P^n}^*$. Dualizing, there is a curve $C$ in $\P^n$ such that all the lines
of $Y$ meet $C$. If $s\ge n$, by dimensional reasons it has to be $s=n$ and
$Y$ will consist of the set of lines meeting $C$. But then $Y$ cannot be
smooth, since the bisecants to $C$ provide singular points of $Y$. Hence,
$s=n-1$ and we can use the classification given in [\ABTcurvafundamental] of
congruences such that all of their lines meet a given curve. However, instead
of checking one by one all the cases in order to find out for which of them
their dual is actually projected from $G(n-2,n+1)$, we will use a more
conceptual method. 

We know from [\ABTcurvafundamental] that $C$ must be a smooth plane curve
--in fact, either a line, or a conic or a smooth plane cubic. In the dual
space, we have then a one-dimensional family of hyperplanes, all of them
containing a fixed $(n-3)$-plane $\Lambda$. If we take a plane $\Pi$ not
meeting $\Lambda$, the intersection of the family of hyperplanes with
$\Pi$ produces a smooth one-dimensional family of lines in $\Pi$. The fact that
the one-dimensional family of hyperplanes is projected from $\P^{n+1}$ is
equivalent to the fact that this curve $C'\subset G(1,\Pi)$ is projected from
$G(1,3)$. But, by Proposition \proydual, this is only possible if $C'$ (and
hence also $C$) is a smooth conic. It also holds that $C'$ is not projected
from $G(1,4)$, so that $Y^*$ is not projected from $G(n-2,n+2)$. Looking again
at the classification in [\ABTcurvafundamental], the only congruences in
$G(1,n)$ for which $C$ is a conic are those of Examples \unodos\ and \dosdos.
We also checked that, for those examples, $Y^*$ can be indeed projected from
$G(n-2,n+1)$.
\qed

\bigskip

\noindent{\sl Remark.} The same classification when $n=3$ was done in [\AS],
where three more cases are found. Even for one of them --the congruence of
bisecants to a twisted cubic-- it holds that
its dual comes not only from $G(1,4)$, but also from $G(1,5)$.

\bigskip

\noindent {\bf Proposition \fibrado-desc}\hfill\break
{\sl Assume $n\ge 4$ and let $Y$ be a congruence in
$G(1,n)$. Then the vector bundle ${\cal S}_{|Y}$ splits if and only if $Y$ is
one of the following:
\item{a)} A congruence contained in some $Y\subset G(1,H)$, for a
hyperplane $H\subset\P^5$.
\item{b)} The congruence of lines joining the points of two disjoint linear
spaces $\Lambda_1$ and $\Lambda_2$ of dimensions $l_1,l_2$ with
$l_1+l_2=n-1$.
\item{c)} The congruence in Example \unodos.}

\bigskip

\noindent{\sl Proof.} Assume ${\cal S}_{|Y}=E_1\oplus E_2$, where each $E_i$
has rank
$e_i$ (hence, $e_1+e_2=n-1$). Since ${\cal S}$ is generated by its sections, so
are
$E_1$ and
$E_2$. This implies that $h^0(Y,E_1)\ge e_1$ and $h^0(Y,E_2)\ge e_2$.
Moreover, if some of the equalities hold, then we have respectively $E_1\cong
{\cal O}_Y^{\oplus e_1}$ or
$E_2\cong {\cal O}_Y^{\oplus e_2}$. But the fact that ${\cal S}_{|Y}$
contains the  trivial line bundle as a summand is equivalent (see Example
\producto) to the fact that
$Y$ is contained in some $G(1,H)$, for a hyperplane $H\subset\P^n$, which is
case a) in the statement. So we can assume $h^0(Y,E_1)\ge e_1+1$ and
$h^0(Y,E_2)\ge e_2+1$.

If $h^0(Y,E_i)=e_i+1$, then it follows that in the dual ${\P^n}^*$ there is a
space of dimension $e_i$ line meeting all the $(n-2)$-planes of $Y^*$ in
dimension $e_i-1$. Dually, there is an $(n-e_i-1)$-plane in
$\P^n$ meeting all the lines of $Y$. Therefore, if both equalities hold, $Y$
must consist of all the lines of $\P^n$ meeting two disjoint linear spaces of
dimensions $(n-e_1-1)$ and $(n-e_1-1)$. This is the congruence of part b) in
the statement, and we checked in Example \producto\ that ${\cal S}_{|Y}$
splits.

So we can assume that, for some $i$, $h^0(Y,E_i)\ge e_i+2$. But
then $h^0(Y,{\cal S}_{|Y})\ge n+2$. Hence $Y$ belongs to one of the
three types of congruences in the statement of Theorem \clas-proy. If $Y$ is
of the first type, i.e. the congruence of lines in a three-dimensional
quadric, as we remarked in Example \Veronese, ${\cal S}_{|Y}$ is
isomorphic to $\Omega_{\P^3}(2)$, hence indecomposable. If $Y$ belongs to
the second type, i.e. it is as in Example \unodos, we checked there that ${\cal
S}_{|Y}$ splits. 

As for the case of the last type of congruences, the one in example \dosdos,
$Y$ is a rational normal scroll of degree $2n-2$, so it spans a linear space of
dimension $3n-4$. This means that, as a projective subvariety under the
Pl\"ucker embedding of $G(1,n)$, it is contained exactly in ${n-2\choose 2}$
linearly independet hyperplanes. The fibers of the scroll consist of the set of
lines contained in a hyperplane $H\subset\P^n$ and passing through a point
$p\in H$. Hence, ${\cal S}_{|Y}$ restricted to a fiber $F\cong\P^{n-2}$ of the
scroll decomposes as ${\cal O}_F\oplus \Omega_F(2)$, and the second summand is
indecomposable. Therefore, if ${\cal S}_{|Y}$ splits as $E_1\oplus E_2$, it
has to be $e_1=1$ and $e_2=n-2$ (or viceversa). Moreover, since $h^0(Y,{\cal
S}_{|Y})=n+2$, either $h^0(Y,E_1)=3$ and $h^0(Y,E_2)=n-1$ or $h^0(Y,E_1)=2$ and
$h^0(Y,E)=n$. Hence, there are two linear subspaces $A_1,A_2\subset\P^n$ of
respective dimensions $r$ and $n-r$ (with $r=1$ or $2$) meeting only at one
point, such that $Y$ is contained in the Schubert varieties of lines meeting
$A_1$ and $A_2$. This implies that $Y$ is contained in
${r\choose 2}+{n-r \choose 2}$ independent hyperplanes, which is a
contradiction, since this number is bigger than ${n-2\choose 2}$ for any value
of $r=1,2$. Therefore, ${\cal S}_{|Y}$ does not split for this last type of
congruences.
\qed

\bigskip

\noindent{\sl Remark.} When $n=3$, there is one more case for which ${\cal
S}_{|Y}$ splits (see [\AS]). Again, this new case is the congruence of
bisecants to a twisted cubic.

\endsection

\vfill\eject 
\centerline{\section References}
\bigskip

\itemitem{\Yo.} {\caps E. Arrondo}, Projections of Grassmannians of lines and
characterization of Veronese varieties, {\it preprint} (1997).

\itemitem{\ABTcurvafundamental.} {\caps E. Arrondo, M. Bertolini and C.
Turrini}, Classification of smooth congruences with a fundamental curve,
pages 43-56 in {\sl Projective geometry with applications} (ed. E. Ballico),
Marcel Dekker, New York, 1994.

\itemitem{\AS.} {\caps E. Arrondo and I. Sols}, {\sl On congruences of lines in
the projective space}, Soc. Math. France (M\'em. {\bf 50}), 1992.

\itemitem{\Rogora.} {\caps E. Rogora}, Varieties with many lines, {\it
Manuscripta Math.} {\bf 82} (1994), 207-226.

\itemitem{\Segre.} {\caps B. Segre}, Sulle $V_n$ aventi pi\`u di $\infty^{n-k}$
$S_k$, I and II, {\it Rend. dell'Acad. Naz. Lincei}, vol. {\bf V} (1948),
193-157 and 217-273.

\itemitem{\Severi.} {\caps F. Severi}, Intorno ai punti doppi impropri di una
superficie generale dello spazio a quattro dimensioni, e a suoi punti tripli
apparenti {\it Rend. Circ. Mat. Palermo II}, {\bf 15} (1901), 377--401.

\itemitem{\Zak.} {\caps F.L Zak}, {\sl Tangents and Secants of Algebraic
Varieties}, Transl. Math. Monographs AMS, vol. {\bf 127}, Providence, RI, 1993

\end